\documentclass{amsart}

\newtheorem{theorem}{Theorem}
\newcommand{\bt}{\begin{theorem}}
\newcommand{\et}{\end{theorem}}
\newtheorem{lemma}{Lemma}
\newcommand{\bl}{\begin{lemma}}
\newcommand{\el}{\end{lemma}}
\newtheorem{corollary}{Corollary}
\newcommand{\bc}{\begin{corollary}}
\newcommand{\ec}{\end{corollary}}
\newcommand{\beq}{\begin{equation}}
\newcommand{\eeq}{\end{equation}}
\newcommand{\benum}{\begin{enumerate}}
\newcommand{\eenum}{\end{enumerate}}
\newcommand{\N}{\ensuremath{ \mathbf N }}
\newcommand{\Z}{\ensuremath{\mathbf Z}}

\newcommand{\mce}{\ensuremath{ \mathcal E}}

\newcommand{\R}{\ensuremath{\mathbf R}}

\newcommand{\Rn}{\ensuremath{ \mathbf{R}^n }}
\newcommand{\Zn}{\ensuremath{ \mathbf{Z}^n }}

\DeclareMathOperator{\conv}{\text{conv}}

\usepackage{amsmath,amssymb,graphicx}

\title{Sums of sets of lattice points and unimodular coverings of polytopes}
\author{Melvyn B. Nathanson}
\address{Department of Mathematics, Lehman College (CUNY), Bronx, NY 10468}
\email{melvyn.nathanson@lehman.cuny.edu}

\subjclass[2010]{Primary  11B13, 11P21, 52A10, 52B20, 52C05.} 
\keywords{Sums of sets of lattice points, lattice polytope, unimodular simplex, 
unimodular covering.}
\thanks{Supported in part by a grant from the PSC-CUNY Research Award Program.}

\date{\today}

\begin{document}

\maketitle

\begin{abstract}
If $P$ is a lattice polytope (that is, the convex 
hull of a finite set of lattice points in $\R^n$),  
then every sum of $h$ lattice points in $P$
is a lattice point in the $h$-fold sumset $hP$.  
However, a lattice point in the $h$-fold sumset $hP$ is not necessarily 
the sum of $h$ lattice points in $P$.
It is proved that if the polytope $P$ is a union of unimodular simplices, 
then every lattice point in the $h$-fold sumset $hP$ is 
the sum of $h$ lattice points in $P$.
\end{abstract}

\section{The addition problem for lattice polytopes and polyhedra} 
The \emph{sumset}, also called the \emph{Minkowski sum}, 
of sets $X_1,\ldots, X_h$ in \Rn\ is the set 
\[
X_1 + \cdots + X_h = \{x_1 + \cdots + x_h:x_i \in X_i \text{ for all } i=1,\ldots, h\}.
\]
If $x_i \in X_i \cap \Z^n$ for $i=1,\ldots, h$, then 
$x_1 + \cdots + x_h \in \Z^n$, and so 
\beq            \label{LatticePolySumset:BasicInclusion}
(X_1 \cap \Z^n) + \cdots + (X_h  \cap \Z^n) 
\subseteq  (X_1 + \cdots + X_h) \cap \Z^n.
\eeq
It is an unsolved problem in additive number theory to describe the $h$-tuples of sets 
$X_1,\ldots, X_h$ in \Rn\ for which equality replaces inclusion 
in the relation~\eqref{LatticePolySumset:BasicInclusion}.  
In particular, this is an unsolved problem in convex geometry,  
even in the important case when the sets $X_1,\ldots, X_h$ are lattice polytopes.

For every positive integer $h$ and every $X \subseteq \R^n$, we define  
the \emph{$h$-fold sumset}
\[
hX = \underbrace{X+\cdots + X}_{\text{$h$ summands}} 
= \{x_1 + \cdots + x_h:x_i \in X\text{ for all } i=1,\ldots, h\}
\]
and, for every positive real number $\lambda$, we define  the \emph{dilation} 
\[
\lambda \ast X = \{ \lambda x: x \in X\}.
\]
If $X$ is convex, then $hX = h \ast X$, that is, the $h$-fold sumset 
equals the dilation by $h$.  

It is also an unsolved problem to determine necessary and sufficient 
conditions for a lattice polytope $P$ to satisfy the equation 
\beq            \label{LatticePolySumset:BasicInclusion-hX}
h(P\cap \Z^n) = (hP) \cap \Z^n  
\eeq
for some integer $h \geq 2$, or for all $h \geq 2$, or for all sufficiently large $h$.  

In this paper we give a simple sufficient condition for a lattice polytope $P$ to satisfy 
equation~\eqref{LatticePolySumset:BasicInclusion-hX} for all positive integers $h$,   
and we show that this sufficient condition implies that for every lattice polytope $P$ 
there is a positive integer $\ell$ such that the lattice polytope $\ell P$ 
satisfies
\beq            \label{LatticePolySumset:BasicInclusion-hlX}
h(\ell P\cap \Z^n) = (h\ell P) \cap \Z^n  
\eeq
all positive integers $h$.

The sufficient condition is that the lattice polytope $P$ have a unimodular cover.  
Both the addition problem for  lattice points in polytopes and unimodular covers and triangulations of polytopes have been extensively investigated.  
For the addition problem, see~\cite{cox-haas-hibi-higa14,haas08,nath-levy16a,oda08}.
For unimodular covers, 
see~\cite{beck-stap10,brun-gube02,brun-gube09,brun-gube-trun97,gube12,karp06,sant-zieg13,thad07}.

\section{Unimodular simplices}
Let $A = \{a_0, a_1,\ldots, a_n\}$ be an affinely independent set in $\R^n$.  
The \emph{$n$-dimensional simplex generated by $A$} is the convex hull of $A$, 
that is, the set  
\begin{align*}
\Delta(A) & = \left\{  \sum_{i=0}^n t_i a_i : t_i \geq 0 \text{ for } i=0,1,\ldots, n \text{ and }\sum_{i=0}^n t_i  = 1 \right\} \\
& = a_0 +  \left\{  \sum_{i=1}^n t_i (a_i - a_0) : t_i \geq 0 \text{ for } i= 1,\ldots, n  \text{ and }\sum_{i=1}^n t_i  \leq 1 \right\}.  
\end{align*}
A \emph{lattice simplex} is an $n$-dimensional simplex $\Delta(A)$,  
where $A$ is a set of $n+1$ affinely intependent lattice points.  

Let $\Delta(A)$ be a  lattice simpliex, and let $\Gamma(A)$ be the subgroup of $\Z^n$ generated by $A-A$.  
The simplex $\Delta(A)$ is \emph{unimodular}\index{unimodular simplex} 
if $\Gamma(A) = \Z^n$.  

For example, the \emph{standard simplex} in $\R^n$ is the lattice simplex 
$\Delta = \Delta\left( \{ 0 \} \cup \mce \right)$,  where  $\mce = \{e_1,\ldots, e_n\}$ 
is the standard basis for $\R^n$.  
This simplex is unimodular.  
In $\R^3$, the simplex generated by the set 
\[
A_1 = \{ 0, e_1, e_2,2e_3\}
\]
 is not unimodular, because 
 \[
 \Gamma(A_1) = \{ (x_1,x_2,x_3) \in \Z^3: x_3 \equiv 0 \pmod{2}\}
 \]
is a subgroup of $\Z^3$ of index 2.  
Note that $\Delta(A_1) \cap \Z^n = A_1 \cup \{e_3\}$.

The simplex generated by the set 
\[
A_2 = \{ 0, e_1, e_2, e_1+e_2+2e_3\}
\]
satisfies $\Delta(A_2) \cap \Z^n = A_2$, 
but $ \Gamma(A_1) =  \Gamma(A_2)$, and so  $\Delta(A_2)$ is not unimodular.

\bl         \label{LatticePolySumset:lemma:unimodular} 
Let $A = \{a_0,a_1,\ldots, a_n\}$ be an affinely independent set in $\Z^n$.  
If the simplex $\Delta(A)$ is unimodular, then
\[
  hA = h\Delta(A) \cap  \Z^n 
\]
every positive integer $h$.  
\el

\begin{proof}
We have $A \subseteq \Delta(A)$ and so $hA \subseteq h\Delta(A)$.
Because $hA \subseteq \Z^n$, it follows that 
$hA \subseteq h\Delta(A) \cap \Z^n$.  

Conversely, let $p \in  h\Delta(A) \cap \Z^n$.  
There exist nonnegative real numbers $t_0, t_1,\ldots, t_n$ such that 
\[
\sum_{i=0}^n t_i = h 
\]
and 
\[
p = \sum_{i=0}^n t_i a_i 
=  \left( h - \sum_{i=1}^n t_i \right) a_0 + \sum_{i=1}^n t_i a_i   
= ha_0 + \sum_{i=1}^n t_i (a_i - a_0) \in \Z^n.
\]
It follows that 
\beq         \label{LatticePolySumset:pha-1}
p - ha_0 = \sum_{i=1}^n t_i (a_i - a_0).
\eeq
The affine independence of the set $A$ implies that the set 
$A-A = \{a_i - a_0:i=1,\ldots, n\}$ 
is an \R-basis for $\R^n$.  
The unimodality of $\Delta(A)$ implies that the set $\{a_i - a_0:i=1,\ldots, n\}$  
is a  $\Z$-basis for the free abelian group $\Z^n$,   
and so there exist integers $w_1,\ldots, w_n$ such that 
\beq         \label{LatticePolySumset:pha-2}
p - ha_0 = \sum_{i=1}^n w_i (a_i - a_0).
\eeq
Comparing equations~\eqref{LatticePolySumset:pha-1} 
and~\eqref{LatticePolySumset:pha-2}, we see that 
\[
t_i = w_i \in \N_0 
\]
for $i= 1,\ldots, n$.  
It follows that 
\[
w_0 = h - \sum_{i=1}^n w_i = h - \sum_{i=1}^n t_i = t_0 \in \N_0.
\]
We have 
\[ 
\sum_{i=0}^n w_i = \sum_{i=0}^n t_i = h    
\]
and so 
\[
p = \sum_{i=0}^n w_i a_i \in hA.
\]
This completes the proof. 
\end{proof}

A similar argument proves that if $\Delta(A)$ is a unimodular simplex in \Rn, 
then $A = \Delta(A) \cap \Z^n$.

\bt
If $P$ is a lattice polytope that is the union of unimodular simplices, then
\[
h\left( P \cap \Z^n \right) =  (hP) \cap \Z^n  
\]
for every positive integer $h$.
Moreover, if $x \in (hP) \cap \Z^n$, then $P  \cap \Z^n$ contains 
an affinely independent set 
$A = \{a_0, a_1, \ldots, a_n\}$ 
such that $x = \sum_{i=0}^n t_ia_i$ with $t_i \in \N_0$ for 
$i=0,1,\ldots, n$ and $\sum_{i=0}^n t_i = h$.  
\et

\begin{proof}
If $p \in (hP) \cap \Z^n$, then $(1/h)p \in P$ and there is a unimodular simplex 
$\Delta(A) \subseteq P$ such that  $(1/h)p \in \Delta(A)$.  
It follows that $p \in h\Delta(A) \cap \Z^n$.
By Lemma~\ref{LatticePolySumset:lemma:unimodular}, 
\[
p \in h\Delta(A) \cap \Z^n = hA \subseteq h(P\cap \Z^n ).
\] 
This completes the proof.  
\end{proof}

\bt
For every lattice polytope $P$ there is a positive integer $\ell$ such that
equation 
\[
h( \ell P\cap \Z^n) = (h \ell P) \cap \Z^n  
\]
holds for every positive integer $h$ that is a multiple of $d$.
\et

\begin{proof}
Every lattice polytope $P$ has a dilation that has a unimodular cover 
(Bruns and Gubeladze~\cite[Chapter 3]{brun-gube09}).  
If $\ell P$ has a unimodular cover, then $\ell P$ 
satisfies equation~\eqref{LatticePolySumset:BasicInclusion-hlX}.  
This completes the proof.  
\end{proof}

 \def\cprime{$'$} \def\cprime{$'$} \def\cprime{$'$} \def\cprime{$'$}
\providecommand{\bysame}{\leavevmode\hbox to3em{\hrulefill}\thinspace}
\providecommand{\MR}{\relax\ifhmode\unskip\space\fi MR }
\providecommand{\MRhref}[2]{%
  \href{http://www.ams.org/mathscinet-getitem?mr=#1}{#2}
}
\providecommand{\href}[2]{#2}


\end{document}